\documentstyle[12pt]{article}

  \newcommand{\const}{\rm const}
  \newcommand{\Var}{\rm Var}
  
  \newcommand{\Law}{\rm Law}

  \newcommand{\Dom}{\rm  Dom}

\begin{document}

\begin{center}

{\bf Ordinary and logarithmical convexity of moment generating function }\\

 \vspace{5mm}

 {\bf M.R.Formica,  E.Ostrovsky and L.Sirota.}

 \end{center}

\vspace{3mm}
 \ Universit\`{a} degli Studi di Napoli Parthenope, via Generale Parisi 13, Palazzo Pacanowsky, 80132,
Napoli, Italy. \\
 \ e-mail: mara.formica@uniparthenope.it \\
\vspace{3mm}

 \ Department of Mathematics and Statistics, Bar-Ilan University, \\
59200, Ramat Gan, Israel. \\
e-mail: eugostrovsky@list.ru\\
\vspace{3mm}

\ Department of Mathematics and Statistics, Bar-Ilan University, \\
59200, Ramat Gan, Israel. \\
e-mail: sirota3@bezeqint.net \\

\vspace{3mm}
\begin{center}
  {\bf Abstract.} \par
\end{center}
\vspace{3mm}
 \hspace{3mm} We establish an ordinary as well as a logarithmical convexity of the  Moment Generating Function
 (MGF) for the centered random variable  and vector (r.v.) satisfying the  Kramer's condition.\par
  \ Our considerations are based on the theory of the so - called Grand Lebesgue Spaces.
\vspace{3mm}
\begin{center}

  {\bf  Key words and phrases.} \par

\end{center}
\vspace{3mm}
 \ Random variable (r.v.), probability, expectation, variance, moment generating function (MGF), ordinary and logarithmical
convexity (OC and LC), Kramer's condition,  derivative, spin glass theory, Gaussian and subgaussian variable and distributions,
parameters, generating and  natural generating function, ordinary and multivariate  Grand Lebesgue Spaces (GLS), norm, tail of
distribution, auxiliary probability space, weight,
slowly varying at infinity function, variable weight function,  Young - Fenchel, or Legendre, transform; examples.\par

\vspace{4mm}

 \section{Statement of problem. Notations. Assumptions.}

\vspace{4mm}

\hspace{3mm} Let $ \ \xi = \xi(\omega) \ $ be a numerical valued {\it centered}: $ \ {\bf E} \xi = 0 \ $ random variable (r.v.)
defined on the suitable sufficiently rich probability space $ \ (\Omega = \{\omega\}, \ {\bf B}, {\bf P}), \ $ having expectation
$ \ {\bf E} \ $ and variance $ \ {\bf \Var}, \ $ satisfying the famous  Kramer's  (Kramer's  - Tchernov's ) condition

\vspace{3mm}

\begin{equation} \label{Kramer cond}
\exists \lambda_0 = \const \in (0,\infty] \ \forall \lambda \in (-\lambda_0, \lambda_0) \Rightarrow \hspace{3mm} {\bf E} \exp(\lambda \ \xi) < \infty.
\end{equation}

\vspace{3mm}

 \ The function

\vspace{3mm}

\begin{equation} \label{MGF}
G(\lambda) = G_{\xi}(\lambda) \stackrel{def}{=} {\bf E} \exp(\lambda \ \xi), \ |\lambda| < \lambda_0
\end{equation}

\vspace{3mm}
is said to be as ordinary as \ {\it  Moment Generating Function } (MGF) for the random variable $ \ \xi. \ $ \par
Introduce also the following functions and set $ \ \Lambda := (-\lambda_0, \lambda_0),  \ $

\vspace{3mm}

\begin{equation} \label{Zete fun}
\Delta(\lambda) = \Delta_{\xi}(\lambda) \stackrel{def}{=} \ln G_{\xi}(\lambda),
\end{equation}

\vspace{3mm}
as well as

\vspace{3mm}

\begin{equation} \label{Phi fun}
\Phi(\lambda) = \Phi_{\xi}(\lambda) \stackrel{def}{=} \lambda^{-1} \ln G_{\xi}(\lambda) = \lambda^{-1} \Delta(\lambda), \ \lambda \ne 0,
\end{equation}

\vspace{3mm}

and $ \ \Phi_{\xi}(0):= 0 = {\bf E}\xi. \ $  \par

\vspace{3mm}

\hspace{3mm} {\bf  Our target in this short report is to establish the convexity of both the introduced functions. } \par

\vspace{3mm}

  \hspace{3mm}  Record: $  \ \Delta_{\xi}(\cdot) \in OC \ $ or equally by notation $ \ \xi \in OC \ $ for the first case
  and correspondingly $ \ \Phi_{\xi}(\cdot) \in LC \ $ or equally $ \ \xi \in LC \ $  in the second one.\par

\vspace{3mm}

 \ The case when the r.v. $ \ \xi \ $ is certain {\it convex} transform  $ \ \xi = F(g) \ $  from a Gaussian centered variable $\ g \ $
 is considered in the recent article \cite{Wei - Kuo Chen}. See also  \cite{Paouris  Valettas},
  \cite{Talagrand 1} - \cite{Talagrand 3}.  \par
 \ Some interest applications if these proposition, especially in the theory of spin glasses, are described in particular in the works
\cite{Charbonneau etc},  \cite{Chatterjee 1} - \cite{Chatterjee 2}, \cite{Comets},
\cite{Contucci and Mingione} -\cite{Guerra 3}, \cite{Jagannath}, \cite{Lehec}, \cite{M�ezard at all}, \cite{Panchenko},
\cite{Parisi 1} - \cite{Parisi 2}, \cite{Sherrington Kirkpatrick} etc. \par

\vspace{4mm}

 \section{Ordinary convexity \ (OC).}

\vspace{4mm}

  \hspace{3mm} {\bf Theorem 2.1.} \ {\sc We conclude under formulated before conditions that the function} $ \ \Delta(\lambda) \ $ {\sc is convex
  on the whole its domain of definition} $ \ \lambda \in \Lambda: \ \xi \in OC. \ $ \par

  \vspace{3mm}

  \ {\bf Proof.} \ We have immediately from  the direct definition  of MGF (\ref{MGF}) for all the values $ \ \lambda \ $ from the set $ \ \Lambda \ $

$$
\Delta'(\lambda) = \frac{{\bf E} \xi \exp(\lambda \xi)}{{\bf E} \exp(\lambda \xi)};
$$

$$
\Delta^{''}(\lambda) = \frac{{\bf E} (\xi^2  \ \exp(\lambda \ \xi)) \ {\bf E} \exp(\lambda \ \xi) - [ {\bf E}\xi \ \exp(\lambda \ \xi) ]^2 }{[{\bf E} \exp(\lambda \ \xi)]^2}.
$$

 \vspace{3mm}

 \ Let us introduce an auxiliary probability space $ \ (\Omega, {\bf B}, {\bf W})  \ $
 grounded on the source {\it measurable} one $ \ (\Omega, {\bf B}), \ $
 with correspondent notions: \ expectation $ \ {\bf H }, \ $  covariation matrix (bilinear operator)  $ \  {\bf Kov}, \ $ grounded on the source
 random  vector $ \ \vec{\xi}, \ $  equipped with a variable positive  weight function  such that for arbitrary  numerical valued random variable
 $ \ \tau  = \tau(\omega) \ $ by definition

 $$
 {\bf H }(\tau) \stackrel{def}{=} \frac{{\bf E} \tau \exp(\lambda \ \xi)}{{\bf E}{\exp(\lambda \xi)}},
 $$
with correspondent variation $ \ {\bf War:} \ $

$$
{\bf War} (\tau) = \frac{{\bf E} \tau^2 \exp(\lambda \ \xi)}{{\bf E}{\exp(\lambda \xi)}} - [{\bf H }(\tau)]^2.
$$
 \ Obviously, both the variables $ \ {\bf H }(\tau), \ {\bf War} (\tau)  \ $ there exists at last in the {\it open} segment
$ \ \lambda \in (\ -\lambda_0, \lambda_0), \ $  for the variables $ \ \tau \ $ for which for instance

$$
|\tau| \le C_0 + C_1 \ |\xi|^k, \ C_0,C_1, \ k \in [0,\infty.)
$$

 \vspace{3mm}
 \ Then the expression for the value $ \ \Delta^{''}(\lambda) \ $ takes the form

$$
\Delta^{''}(\lambda) = {\bf War} (\xi),
$$
which is  finite and non - negative for all the admissible values  $ \ \lambda \in \Lambda; \ $ \par
Q.E.D.\par

\vspace{4mm}

 \section{Logarithmical convexity (LC). }

\vspace{4mm}

 \begin{center}

{\it Grand Lebesgue Spaces (GLS). }

 \end{center}

 \ Let $ \ \phi = \phi(\lambda), \ \lambda \in \Lambda = (-\lambda_0, \ \lambda_0), \ \lambda_0 = \const \in (0, \infty] \ $  be certain
 numerical valued even strong convex which takes only positive values for positive arguments twice continuous differentiable
function, briefly: Young-Orlicz function, and such that

$$
\phi(0) = \phi'(0) = 0, \  \exists \ \phi''(0) \in (0,\infty).
$$

 \hspace{3mm}  For instance: $ \ \phi(\lambda) = 0.5\lambda^2, \ \lambda_0  = \infty; \ $ the so-called {\it subgaussian case}.\par

 \ We denote the set of all these Young-Orlicz function as $ \ F; \ F = F_{\Lambda} = \{\phi(\cdot)\} = \{\phi[\Lambda](\cdot) \ \}. \ $ \par

 \ We say by definition that the {\it centered:} \ $ \ {\bf E} \xi \ $ random variable (r.v) $ \ \xi  = \xi(\omega) \ $ belongs
to the space $ \ B(\phi), \ \phi(\cdot) \in F \ $ iff there exists certain non-negative {\it constant} $ \ \rho \ge 0 \ $ such that

$$
\forall \lambda \in (-\lambda_0, \lambda_0) \ \Rightarrow \max_{\pm}
E \exp( \ \pm \ \lambda \ \xi) \le \exp[\phi(\lambda \ \rho)].
$$

\vspace{3mm}

 \ The minimal {\it non-negative} value $ \ \rho \ $  satisfying the last estimate for all the admissible values
$ \ \lambda \in ( - \lambda_0, \ \lambda_0), \ $  is named as a $ \ B(\phi) \ $ norm of the random variable $ \ \xi \ $, write

$$
||\xi||B(\phi) \stackrel{def}{=} \inf \{\rho, \ \rho > 0: \ \forall \lambda \in (-\lambda_0, \ \lambda_0) \ \Rightarrow
\max_{\pm} {\bf E} \exp(\ \pm \lambda \ \xi ) \le  \exp(\phi(\lambda \ \rho)) \ \}.
$$

\vspace{3mm}

 \ The set $ \ B(\phi) \ $ relative the introduced norm $ \ ||\cdot||B(\phi) \ $ and equipped with ordinary algebraic operations forms
 a (complete) Banach functional  rearrangement invariant space. The function $ \ \phi(\cdot) \ $ is named as usually as {\it generating function}
 for this space. These spaces are named as usually {\bf  Grand Lebesgue Spaces \ (GLS)}.  \par

 \vspace{3mm}

 \ These spaces are very convenient for the investigation of the r.v. having an
exponential decreasing tail of distribution, for instance, for investigation of the limit
theorem, the exponential bounds of distribution for sums of random variables, non - asymptotical properties, problem of continuous
and weak compactness of random fields, study of Central Limit Theorem in the Banach space etc. \par

\vspace{3mm}

 \hspace{3mm} See about this definitions the seminal works \cite{Ermakov},
\cite{Fiorenza1} - \cite{Fiorenza4}, \cite{formicagiovamjom2015},  \cite{Kozatchenko  Ostrovsky 1},
\cite{Kozatchenko  Ostrovsky 2}, \cite{Ostrovsky 1},  \cite{Ostrovsky 2}. \par

\vspace{3mm} A very important example: the - called subgaussian random variables. Let $ \ \Lambda = R \ $ and
$ \ \phi(\lambda) = \phi_2(\lambda) \stackrel{def}{=} 0.5 \ \lambda^2; \ $ see e.g.
\cite{Buldygin Kozachenko}, \cite{Kahane 1} - \cite{Kahane 2}. \par

\vspace{3mm}

 \ {\bf  Definition 3.1. } \ {\it  A natural generating function. \ \ } Let the random variable $ \ \eta \ $ be some {\it centered} non - zero r.v. satisfying
 the famous Kramer's  (Kramer's \ - \ Tchernov's)  condition:

$$
\exists C \in (0,\infty) \ \Rightarrow {\bf P}(|\eta| > t) \le \exp(- C t), \ t \ge 0;
$$

 or equally

$$
\exists C_0 > 0 \ \ \forall \mu: \ |\mu| < C_0 \ \Rightarrow  {\bf E}  \exp(\mu \ \eta) < \infty.
$$

\vspace{3mm}

 \ The following generating function

$$
\nu(\lambda) = \nu[\eta](\lambda) \stackrel{def}{=} \ln \max({\bf E}\exp(\lambda \ \eta), \ {\bf E}\exp(-\lambda \ \eta) \ )
$$
is named as usually as a {\it natural generating function } for the r.v. $ \ \eta. \ $ \par
 \ Obviously, this function is finite (and convex) inside some non - trivial neighborhood of origin. \par

\vspace{3mm}

 \ Of course, $ \  \eta \in B(\nu[\eta]) \ $ and $ \ ||\eta|| B(\nu[\eta]) = 1. \ $ \par

\vspace{3mm}

 \hspace{3mm} {\it An example.}  \ Let   the {\it symmetrical} distributed r.v. $ \ \zeta = \zeta[m,\gamma,L] \ $ be such that
$ \ \nu[\zeta](\lambda) = \phi[m,\gamma,L](\lambda), \ $ where by definition

\vspace{3mm}

\begin{equation} \label{mult example}
 \phi[m,\gamma,L](\lambda) \asymp  C_1 \ |\lambda|^m \ \ln^{\gamma}(|\lambda|) \ L(\ln |\lambda|), \  |\lambda|  > Z = \const > e,
\end{equation}

\vspace{3mm}

 and of course

$$
\phi[m,\gamma,L](\lambda) \asymp C_0 \ \lambda^2, \  |\lambda| \le Z, \  C_0, C_1 \in (0,\infty).
$$

\vspace{3mm}

\ Here $ \ m = \const > 0, \ \gamma = \const \in R; \ L(\cdot) \ $ be positive continuous slowly varying at infinity function.\par
\ Many examples of these functions which are MGF functions for certain random variables
are represented in   \cite{Kozatchenko  Ostrovsky 2}, \cite{Ostrovsky 1},  \cite{Ostrovsky 2}.\par

\vspace{3mm}

  \ We will write for these r. v.  $ \  \Law(\zeta) = G(m,\gamma,L); \ $ but in the case when $ \ L(\lambda)= 1, \ |\lambda|> Z \ $ we denote
for the sake of simplicity $ \  \Law(\zeta) = G(m,\gamma). \ $\par

\vspace{3mm}

\  See the detail investigations of these spaces in the many  works \cite{Ermakov}, \cite{Fiorenza1}, \cite{Fiorenza2},
 \cite{Fiorenza3}, \cite{Fiorenza4}, \cite{formicagiovamjom2015}, \cite{Kozatchenko  Ostrovsky 1}, \cite{Kozatchenko  Ostrovsky 2}.\par
\ In particular, there are  studied  therein the  {\it exact}  interrelations between the tail behavior of the random variable, its moment behavior and
belonging of these r.v. to certain $ \ B(\phi) \ $ space.  For instance, for the  {\it centered}  r.v. $ \  \delta \ $ there holds the following implications

$$
\ln {\bf P} (|\delta| > x)  \sim - c_1 x^d, \ d = \const > 1, \ x \ge e, \ x \to \infty,\ \Longleftrightarrow
$$

$$
\ln {\bf E}\exp ( \pm \lambda \ \delta) \sim c_2 |\lambda|^{d/(d-1)}, \ |\lambda| \to \infty, \Longleftrightarrow
$$

$$
\frac{||\delta||_p}{p^{1/d} } \sim c_3, \ p \to \infty;  \ ||\delta||_p \stackrel{def}{=} ( {\bf E}|\delta|^p )^{1/p}, \ p \ge 1.
$$

\vspace{3mm}

 \ Here $ \ c_1,c_2,c_3 = c_1,c_2,c_3(d) = \const \in (0,\infty). \ $\par

\ The case $ \ d = 2 \ $ correspondent to the well known subgaussian random variable (distribution).\par

\vspace{4mm}

 \hspace{3mm} {\sc  Proposition 3.1.} \par

 \hspace{3mm} To be more generally, suppose that the non - zero centered r.v. $ \  \beta \ $  belongs to some $ \ B(\phi)  \ $  space; $ \ \phi \in F. \ $
 Let for definiteness $ \ ||\beta||B(\phi) = 1. \ $  Then

\begin{equation} \label{YFT}
\max  \left[ \ {\bf P}( \beta > x ), \ {\bf P}( - \beta > x )  \ \right]  \le \exp \left(-\phi^*(x)  \ \right), \ x > 0.
\end{equation}

\vspace{3mm}

 \ Here  $ \ \phi^*(\cdot)  \ $  denotes the so - called Young - Fenchel, or Legendre, transform for the function $ \ \phi: \ $

$$
\phi^*(x)  \stackrel{def}{=} \sup_{y \ge 0} (xy - \phi(y)).
$$

\vspace{3mm}

\ Herewith the inverse concludion to (\ref{YFT}) also holds true, of course, up to multiplicative finite positive constant.  Indeed, let
(\ref{YFT}) take place, then under appropriate restrictions  \cite{Kozatchenko  Ostrovsky 1}, \cite{Kozatchenko  Ostrovsky 2}

\begin{equation} \label{MGF estim}
{\bf E } \exp(\lambda \beta) \le \exp \phi( C \ \lambda ), \ C = \const > 0
\end{equation}
for some non - trivial centered neighborhood for the values  of parameter $ \ \lambda: \ $

$$
\exists \ \Lambda_1 = \const \in (0,\infty] \ \Rightarrow |\lambda| \in [0,\Lambda_1),
$$
see yet \cite{Bagd Ostr}. Equally: \ $ \ ||\beta|| B(\phi) < \infty. \ $ \par

\vspace{3mm}

\ {\sc Remark 3.0. \ } This assertion remains true still in the multidimensional case. Indeed, the classical Youg - Fenhel (or Legendre)
transform has  here a form

$$
\phi^*(\lambda) \stackrel{def}{=} \sup_{x \in R^n}[ (x,\lambda) - \phi(x) ],
$$
where $ \ x =\vec{x} \in R^n, \ \lambda = \vec{\lambda} \in R^n,  \ $  and as ordinary $ \  (x, \lambda)  \ $ denotes the inner (scalar)
product. \par

\vspace{4mm}

 \hspace{3mm} {\sc Proposition 3.2. \  } \par

 \vspace{3mm}

{\sc Let }  $ \  \Law(\zeta) = G(m,\gamma). \ $  { \sc Then  the distribution of this r.v. is logarithmical convex:} $ \ \zeta \in LC \ $
{\sc  if and only if  either } $ \ m > 1 \ $ {\sc or} $ \ m = 1, \ \gamma \ge 0. \ $  \par

\vspace{3mm}

\ Indeed, the convexity of this function in the domain $ \ [-Z,Z] \ $  is obvious. Further,
it is easily to verify that the function $ \ G(m, \gamma)(\lambda)/\lambda \ $ is convex in the rest of semi - axis  $ \ [Z,\infty) \ $
only under formulated above conditions. \par

 \vspace{4mm}

  \section{Multivariate case.}

 \vspace{4mm}

 \hspace{3mm} It is interest in our opinion to generalize obtained results into the multidimensional case.  Namely, let
 $ \ \xi  = \vec{\xi} \ = \{ \ \xi_1, \xi_2, \ \ldots, \xi_l \ \}, \ l = 2,3, \ldots \  $ be certain {\it centered} random vector satisfying
 a so - called multidimensional analog of  the Kramer's condition

 $$
  \exists \delta_0 > 0, \ \forall \lambda = \vec{\lambda}: \ |\lambda| < \delta_0   \Rightarrow
 $$

$$
Q(\lambda) = Q_{\xi}(\lambda)  \stackrel{def}{=} {\bf E} \exp(\vec{\xi}, \vec{\lambda}) < \infty.
$$

 \ Here as ordinary  $ \ \lambda = \vec{\lambda} = \{ \lambda_1, \lambda_2, \ldots, \lambda_l \},  $

 $$
 |\lambda| = \sqrt{ \sum_{i=1}^l \lambda_i^2 }, \hspace{3mm} (\lambda,\xi) = \sum_{i=1}^l \lambda_i \ \xi_i.
 $$

\vspace{3mm}

 \ The function $ \ Q(\lambda) = Q_{\xi}(\lambda) \ $ is named as before as Moment Generating Function (MGF) for the random vector $ \ \vec{\xi}. \ $ \par
 \ The convexity of the function  $ \ \tilde{Q}(\lambda) \stackrel{def}{=} \ln \left[ \ Q_{\xi}(\lambda) \ \right], \ \lambda \in \Dom \{Q\} $ may be grounded quite alike  to the
 one - dimensional case $ \ l = 1. \ $\par
 Consider now the following analog for its "logarithmic  derivative"  \  (LD), where by definition

$$
V(\lambda) = V_{\xi}(\lambda) \stackrel{def}{=} \ln Q_{\xi}(\lambda)/|\lambda|, \ \vec{\lambda} \ne 0; \ V(0) \stackrel{def}{=} 0;
$$
and will write $ \  \vec{\xi} \in LD, \ $ iff this function $ \ V(\lambda) \ $ is convex in the whole of its domain of finiteness.\par
 \ Let us consider the case when the centered random vector  $ \ \xi = \vec{\xi} \ $ has a following type of MGF

$$
\ln {\bf E} \exp(\lambda,\xi) \asymp |\lambda|^m\ \ln^{\gamma}(|\lambda|), \ |\lambda| > e.
$$
  Here as above $ \ m = \const > 0, \ \gamma = \const \ \in R. \ $ We conclude that $ \  \vec{\xi} \in LD, \ $ iff
 $ \ m > 1 \ $ or $ \ m = 1,\ \gamma > 0. \ $ \par

\vspace{3mm}

 \ {\it A subexample:  multivariate subgaussian random vector.}  Let $ \  B  \ $  be symmetrical  strictly positive definite  deterministic matrix
 of an order $ \  l  \times l, \ l = 2,3, \ldots. \ $ The  centered random vector $ \ \zeta = \vec{\zeta} \ $  is said to be  by definition
subgaussian relative this matrix $ \ B, \ $   iff

$$
\ln {\bf E} \exp(\lambda,\zeta) \asymp 0.5 \ (B\lambda,\lambda), \ \lambda \in R^l.
$$

 \vspace{3mm}

\ As long as here $ \ m = 2, \ $ we conclude  that  $ \ \zeta \in LD. \ $ \par

\vspace{6mm}

\vspace{0.5cm} \emph{Acknowledgement.} {\footnotesize The first
author has been partially supported by the Gruppo Nazionale per
l'Analisi Matematica, la Probabilit\`a e le loro Applicazioni
(GNAMPA) of the Istituto Nazionale di Alta Matematica (INdAM) and by
Universit\`a degli Studi di Napoli Parthenope through the project
\lq\lq sostegno alla Ricerca individuale\rq\rq .\par


\vspace{6mm}


 \begin{thebibliography}{44}

\bibitem{Bagd Ostr}
{\bf D. R. Bagdasarov,  \ E.I. Ostrovsky.} {\it Exponential Confidence Intervals in the Nonparametric Density Estimation.} \
Theory of Probability and its Applications, 2012, (1997), SCI 4, pp. 47 - 63.
Obninski Instute for Nuclear  Power Energy.

\bibitem{Bobkov 1}
 {\bf S. Bobkov.} {\it Spectral gap and concentration for some spherically symmetric probability measures.} In: \
 Geometric Aspects of Functional Analysis: Israel Seminar 2001-2002, pages 37 \ - \ 43. Springer, 2003.

\bibitem{Bobkov 2}
{\bf S. Bobkov, G. Chistyakov, and F. G\"otze.} {\it Strongly subgaussian probability distributions.} Preprint, 2022.


\bibitem{Buldygin Kozachenko}
{\bf V.V. Buldygin and Yu.V. Kozachenko.} {\it Subgaussian random variables.}
Ukrainian Math. J.,  {\bf 32,} \ (1980), 483 \ - \ 489


\bibitem{Charbonneau etc}
{\bf P. Charbonneau, E. Marinari, G. Parisi, F. Ricci-Tersenghi, G. Sicuro, F. Zamponi, and M. M�ezard.}
{\it Spin Glass Theory and Far Beyond: Replica Symmetry Breaking after 40 Years.} World Scientific, 2023.


\bibitem{Chatterjee 1}
{\bf S. Chatterjee.} {\it Chaos, concentration, and multiple valleys.}  arXiv preprint; \ arXiv:0810.4221, 2008.

\bibitem{Chatterjee 2}
{\bf S. Chatterjee.}  {\it Disorder chaos and multiple valleys in spin glasses.} arXiv preprint; \ arXiv: 0907.3381, 2009.

\bibitem{Comets}
{\bf F. Comets.}  {\it Directed polymers in random environments.} Lecture Notes in Mathematics, 2017.

\bibitem{Contucci and Mingione}
{\bf P. Contucci and E. Mingione.}  {\it A multi-scale spin-glass mean-field model.} Communications in Mathematical Physics, 368: 1323 \ - \ 1344, 2019.

\bibitem{Ermakov}
{\bf V. Ermakov, and E. I. Ostrovsky.} {\it Continuity Conditions,  Exponential Estimates, and the Central Limit Theorem for Random Fields.}
Moscow, VINITY, (1986), (in Russian).


\bibitem{Fiorenza1}
{\bf Fiorenza A, and Karadzhov G.E.} {\it Grand and small Lebesgue spaces and their analogs.} Journal
for Analysis and its Applications 2004; 23 (4) : 657 \ - \ 681.

\bibitem{Fiorenza2}
 {\bf Fiorenza A.} {\it Duality and reflexity in grand Lebesgue spaces.} Collect. Math. 2000; 51 ( 2) : 131 \ - \ 148.


 \bibitem{Fiorenza3}
{\bf A.~Fiorenza, M. R.~Formica, A.~Gogatishvili, T.~Kopaliani} and  {\bf J.~M.Rakotoson.}
 {\it Characterization of interpolation between grand, small or classical Lebesgue spaces}. \\
arXiv:1812.04295v1 [math.FA] 11 Dec 2018

\bibitem{Fiorenza4}
{\bf Fiorenza A., and Karadzhov G.E.} {\it Grand and small Lebesgue spaces
and their analogs.} Consiglio Nationale Delle Ricerche, Instituto per le Applicazioni
del Calcoto Mauro Picine, Sezione di Napoli, Rapporto tecnico n. 272/03, (2005).

\bibitem{formicagiovamjom2015}
{\bf M.~R. Formica} and {\bf R.~Giova.} {\it Boyd indices in
generalized grand Lebesgue spaces and applications}. Mediterr. J.
Math. \textbf{12} (2015), no.~3, 987--995.


\bibitem{Guerra 1}
{\bf  F. Guerra.}  {\it Broken replica symmetry bounds in the mean field spin glass model.} Communications in
Mathematical Physics, 233: 1 \ - \ 12, 2003.

\bibitem{Guerra 2}
{\bf F. Guerra.} {\it The emergence of the order parameter in the interpolating replica trick for disordered
statistical mechanics systems.} In: Lake Como School of Advanced Studies Complexity and Emergence:
Ideas, Methods, with special attention to Economics and Finance, pages 63 \ - \ 86. Springer, 2018.


\bibitem{Guerra 3}
{\bf F. Guerra and F. L. Toninelli.} {\it The thermodynamic limit in mean field spin glass models.}
 Communications in Mathematical Physics, 230: 71 \ - \ 79, 2002.

\bibitem{Jagannath}
{\bf A. Jagannath.}  {\it Approximate ultrametricity for random measures and applications to spin glasses.}
 Communications on Pure and Applied Mathematics, 70(4): 611 \ - \ 664, 2017.



\bibitem{Kahane 1}
{\bf J.P. Kahane.} {\it Proprietes locales des fonctions l`a series de Fourier al'eatoires.} Stud. Math. 19 (1960), 1 \ - \ 25.

\bibitem{Kahane 2}
{\bf J.P. Kahane.}  {\it Some random series of functions.} 2nd ed. Cambridge
University Press, London, 1985.

\bibitem{Kozatchenko  Ostrovsky 1}
{\bf Kozatchenko Yu. V., Ostrovsky E.I.} {\it The Banach Spaces of random Variables of subgaussian type. }
Theory of Probab. and Math. Stat. (in Russian). Kiev, KSU, {\bf 32,} \ 43 \ - \ 57, (1985).

\bibitem{Kozatchenko  Ostrovsky 2}
{\bf Kozachenko Yu.V., Ostrovsky E., Sirota L.} {\it Relations between exponential tails, moments and
moment generating functions for random variables and vectors.} \\
arXiv:1701.01901v1 [math.FA] 8 Jan 2017


\bibitem{Lehec}
{\bf J. Lehec.}  {\it Short probabilistic proof of the Brascamp-Lieb and Barthe theorems.} Canadian Mathematical
Bulletin, 57(3): 585 \ - \ 597, 2014.

\bibitem{M�ezard at all}
{\bf M. M�ezard, G. Parisi, and M. A. Virasoro.}  {\it Spin glass theory and beyond: An Introduction to the Replica
Method and Its Applications,} volume 9. World Scientific Publishing Company, 1987.

\bibitem{Ostrovsky 1}
{\bf Ostrovsky E.I. } (1999). {\it Exponential estimations for Random Fields and its
applications.} (in Russian). Moscow \ - \ Obninsk, OINPE.

\bibitem{Ostrovsky 2}
{\bf Ostrovsky E. and Sirota L.} {\it Vector rearrangement invariant Banach spaces
of random variables with exponential decreasing tails of distributions.} \\
 arXiv:1510.04182v1 [math.PR] 14 Oct 2015

\bibitem{Panchenko}
{\bf Dmitry Panchenko.}  {\it The Parisi formula for mixed p-spin models.} The Annals of Probability, pages
946 \ - \ 958, 2014.


\bibitem{Paouris  Valettas}
{\bf G. Paouris and P. Valettas.} {\it A Gaussian small deviation inequality for convex functions.} The Annals of
Probability, 46(3): 1441 \ - \ 1454, 2018.

\bibitem{Parisi 1}
{\bf G. Parisi.} {\it The order parameter for spin glasses: a function on the interval 0-1.} Journal of Physics A:
Mathematical and General, 13(3): 1101 - 1111, 1980.

\bibitem{Parisi 2}
{\bf G. Parisi.}  {\it A sequence of approximated solutions to the SK model for spin glasses.} Journal of Physics
A: Mathematical and General, 13(4): L115, 1980.

\bibitem{Sherrington Kirkpatrick}
{\bf D. Sherrington and S. Kirkpatrick.}  {\it Solvable model of a spin-glass.} Phys. Rev. Lett., 35: 1792 \ - \ 1796, Dec. 1975.

\bibitem{Talagrand 1}
{\bf M. Talagrand.} {\it The Parisi formula.} Annals of Mathematics, pages 221 \ - \ 263, 2006.

\bibitem{Talagrand 2}
{\bf M. Talagrand.}  {\it  Large Deviations, Guerra�s and ASS Schemes, and the Parisi Hypothesis.} Journal of
Statistical Physics, 126, 2007.


\bibitem{Talagrand 3}
{\bf M. Talagrand.} {\it Mean field models for spin glasses.} Volume I: Basic examples. Volume 54. Springer
Science  and Business Media, 2010.



\bibitem{Wei - Kuo Chen}
{\bf Wei-Kuo Chen.} \ {\it A Gaussian Convexity for Logarithmic Moment Generating Functions.} \\
arXiv:2311.0835 v1 [math. PR] 14 Nov 2023

\end{thebibliography}
\end{document}